\def\ZZ         {{\mathbb Z}}

\def\CC         {{\mathbb C}}
\def\QQ         {{\mathbb Q}}
\def\PP         {{\mathbb P}}

\def\A          {{\cal A}}

\def\J          {{\cal J}}

\def\O          {{\cal O}}

\def\S           {{\cal S}}

\def\Char       {{\rm Char}}

\def\cal        {\mathcal}

\documentclass{amsart}
\usepackage{amssymb}
\usepackage{verbatim}
\usepackage{eucal}
\usepackage{mathrsfs}
\usepackage{graphicx}
\usepackage{psfrag}

\newtheorem{theorem}{Theorem}[section]

\newtheorem{prop}[theorem]{Proposition}

\theoremstyle{definition}
\newtheorem{dfn}[theorem]{Definition}

\newtheorem{example}[theorem]{Example}

\theoremstyle{remark}
\newtheorem{remark}[theorem]{Remark}

\begin{document}

\title[On combinatorial Invariance and Catalan equation]
{On combinatorial invariance of the cohomology of
Milnor fiber of arrangements and 
Catalan equation over function fields.}

\author{A.Libgober}

\address{
Department of Mathematics\\
University of Illinois\\
Chicago, IL 60607 and 
Department of Mathematics \\
University of Miami \\
Coral Gables, Fl.}
\email{libgober@math.uic.edu}

\thanks{The author was partially supported by NSF grant.}

\begin{abstract} 
We discuss combinatorial invariance
of the betti numbers of the Milnor fiber 
for arrangements of lines 
with points of multiplicity two and three
and describe a link between this problem 
and enumeration of solutions of the Catalan 
equation over function field in the case when its coefficients 
are products of linear forms and the equation 
defines an elliptic curve.  
\end{abstract}

\maketitle

\section{Introduction}

An interesting problem in the theory of arrangements 
of hyperplanes over the field of complex numbers 
is to determine which topological invariants of the 
complement depend only on the combinatorics
of the arrangements i.e. the poset 
formed by intersections 
of the hyperplanes of the arrangement. A theorem, 
going back to Arnold, Orlik and Solomon (cf. \cite{OS}), asserts that the
cohomology of the complement is a combinatorial invariant.
On the other hand, it is known that the homotopy 
type (cf. \cite{Ryb}) and even the fundamental group (cf. \cite{cogoryb}) 
cannot be determined from the combinatorics of the 
arrangement since there exist pairs of arrangements 
with the same combinatorics but different such invariants.

Given an arrangement $\A$ of hyperplanes in $\PP^{n+1}$ in which 
the equations of hyperplanes are $L_i=0$ where $L_i$ are 
linear forms, one can consider the affine hypersurface 
$\Pi L_i=1$ which is called the Milnor fiber of arrangement.
It is known that for the Milnor fiber $M_{\A}$ the homology 
group $H_1(M_{\A},\ZZ)$ depends only on the fundamental group
of the complement to $\A$ 
and more precisely it is the abelianization 
of the subgroup of $\pi_1(\PP^2-\A)$ which the kernel of the 
surjection $lk: \pi_1(\PP^2-\A) \rightarrow \ZZ_{{\rm Card} \A}$ given 
by the total linking number with $\A$ cf. (\cite{charvar}).
One can ask about 
the combinatorial invariance of the cohomology of $M_{\A}$ and in 
particular $rk H^1(M_{\A},\QQ)$
called (the first) Milnor number of the arrangement
(for arrangements of lines, this is equivalent to the combinatorial 
invariance of remaining cohomology groups). 
To this end we show the following:

\begin{theorem}\label{invariance} If two arrangements of lines 
$\A_1$ and $\A_2$ with the multiple points
of multiplicity three and two only  
are combinatorially equivalent and if the monodromy 
on the Milnor fiber of $\A_1$ has an eigenvalue different from one
then so does the monodromy on the Milnor fiber the arrangement $\A_2$.
\end{theorem} 

This theorem will be derived from the following more precise
statement relating the cohomology of Milnor fiber 
and other data of arrangement.
We say that an arrangement $\A$ is composed of a pencil
$L$ of plane curves (cf.\cite{withsergey},\cite{falksergey})
if $L=\{\lambda F+\mu G \vert deg F=deg G\}$ contains  
three curves $\lambda F+\mu G$ 
with equations which are products of linear forms 
such that $\A$ is the zero set of the forms defining these 
curves. A pencil is reduced if 
equations of the above reducible curved are reduced.
(see section \ref{prelim} for other definitions needed for this theorem). 

\begin{theorem}\label{corrected} 
Let $\A$ be an arrangement of lines with double 
and triple points only.
If the monodromy acting on  
one-dimensional cohomology of the Milnor number of $\A$ 
has an eigenvalue different from one then the 
arrangement $\A$ is composed of a reduced pencil.
The characteristic variety of the arrangement 
$\A$ has a component containing the identity
$Char \pi_1(\PP^2-\A)$.
Vice versa, if $\A$ is composed of a reduced pencil then  
the monodromy has eigenvalue $exp({{2\pi i} \over 3})$.
\end{theorem}

The argument for the proof of above results 
relies on the relation between the Alexander polynomials 
of plane curves and Mordell-Weil groups 
considered recently in \cite{jose} and study 
of the Catalan equation over function field.
Catalan equation, i.e. $X^p-Y^q=1$ or more generally 
\begin{equation}\label{catalangeneral}
A_1X^p+A_2Y^q+A_3Z^r=0
\end{equation}
has a long history and its properties 
of course depend on the ground field. Over function field,
i.e. the case when $A_i \in \CC({\cal X})$ where 
$\cal X$ is a projective manifold 
the equation (\ref{catalangeneral}) was considered by 
J.Silverman. A finiteness results was shown in \cite{silverman} 
in the case the equation (\ref{catalangeneral}) defines over 
the function field a curve of genus greater than 1.
Elementary treatment of Fermat 
equation (i.e. $A_i \in \CC, p=q=r$) is given in \cite{green}. 
Finding solution in the case 
of genus zero is trivial. Present paper considers the
case $p=q=r=3$ and $A_i$ are product of linear forms in 3 variables.
The answer to question about enumeration of solutions
depends on the geometry of the curve $A_1 \cdot A_2\cdot A_3=0$ 
in $\PP^2$ and we describe the cases when the number of solutions is infinite 
(cf. theorems \ref{catalan} and \ref{catalanonevar}).

The paper organized as follows. In the next section we recall terminology
from the theory of characteristic varieties of arrangements and 
related background information. In section \ref{catalansection} 
we prove results on existence and enumeration of solutions 
to Catalan equation. Proof of theorems  presented 
in the last section. The paper ends with remark on 
properties of Mordell Weil groups which would imply combinatorial 
invariance of the cohomology of Milnor fiber.
The author wants to thank S.Kaliman and L.Katzarkov for 
hospitality during visit to University of Miami where 
this work was written as well as S.Yuzvinsky for a useful comments.
Finally, I grateful to A.Dimca for pointing 
out example \ref{dualflex} which helped to correct error in 
original statement of theorem \ref{corrected} and to J.I.Cogolludo-Agustin 
for numerous examples of Catalan equations 
which helped to correct theorem \ref{catalan}.

\section{Preliminaries}\label{prelim}

As was already mentioned, Milnor fiber of an arrangement $\A$ of hyperplanes 
in $\PP^{n+1}$ is affine hypersurface $M_{\A}$ in $\CC^{n+2}$ given by 
equation $\Pi L_i=1$ where $L_i$ are 
linear forms defining the hyperplanes of the 
arrangement. The standard map $\CC^{n+2} -{\cal O} \rightarrow \PP^{n+1}$
($\O$ is the origin) 
restricted on $M_{\A}$ is a cyclic cover 
\begin{equation}\label{cover}
       M_{\A} \rightarrow \PP^{n+1}-\A
\end{equation}
of degree equal to 
the number of hyperplanes $r=Card \A$ in $\A$. This cyclic cover 
corresponds to the homomorphism $\pi_1(\PP^{n+1}- \A) \rightarrow 
\ZZ_{r}$ onto the cyclic group of order $r={\rm Card} \A$ 
given by the total modulo $r$ linking number with the union of hyperplanes in 
$\A$. We are interested in $dim H_1(M_{\A},\CC)$ and also in 
the eigenspaces of the generator of the group $\ZZ_r$ of 
deck transformations acting on $H_1(M_{\A},\CC)$. 
By Lefschetz theorem applied 
to a generic plane $H$ in $\PP^{n+1}$ this group is isomorphic to the 
corresponding group of line arrangement formed by intersections
of $H$ with the hyperplanes of $\A$. So from now on we shall assume that 
$n=1$. Since the main results in this note 
deal with the case of arrangements of lines 
with points of multiplicity 2 and 3 from now on we shall restrict ourselves 
 to the case of such  multiplicities as well. 

By \cite{duke}, \cite{arcata} 
(cf. also \cite{charvar}) 
the rank of the cyclic cover (\ref{cover}) can be found in terms the 
Alexander polynomial of $\PP^2-\A$. More precisely, the Alexander polynomial 
of $\PP^2-\A$  
is equal to 
\begin{equation}\label{alex}
(t-1)^{{\rm Card} \A-2}(t^2+t+1)^{s}
\end{equation}
for $s=dim H^1(\PP^2,{\cal J}_{\S}({{2 {\rm Card \A}} \over 3}-3))$
where $\J_{\S} \subset \O_{\PP^2}$ is the ideal sheaf having stalk different 
from $\O_P$ iff $P$ is a triple point of the arrangement and 
equal to maximal ideal at those $P$. Moreover if $r$ is divisible by $3$ then 
\begin{equation}
{\rm rk}H_1(M_{\A})_{\omega}={\rm rk}H_1(M_{\A})_{\omega^2}=s, \ \ {\rm rk}H_1(M_{\A})_1=r-2  
\end{equation}
where $H_1(M_{\A})_{\zeta}$ 
is the eigenspace of the deck transformation of the 
covering (\ref{cover}) with eigenvalue $\zeta$ and $\omega=exp({{2 \pi \sqrt{-1}} \over 3})$.
If $r$ is not divisible by 3 then the dimensions of eigenspaces with 
eigenvalue $\omega, \omega^2$ are equal to zero. Hence from now on we shall 
assume that $r$ is divisible by $3$.

Together with unbranched cover $M_{\A}$ of degree $r$ one can consider 
a non-singular model of {\it 3-fold} branched 
cyclic covering $\bar M_{\A}$ 
of $\PP^2$ which is a compactification of the quotient of $M_{\A}$ 
by the action of the subgroup $\ZZ_{r \over 3} \subset \ZZ_r$.
The integer ${\rm rk} H_1(\bar M_{\A},\CC)$ is birational 
invariant and hence depends
only on $M_{\A}$. 
Using results of \cite{duke}, one has:  
\begin{equation}\label{dimeigenspace}
{1 \over 2}{\rm rk}H_1(\bar M_{\A})={\rm rk}H_1(\bar M_{\A})_{\omega}=
{\rm rk}H_1(\bar M_{\A})_{\omega^2}=s
\end{equation}

Recall the following definition from \cite{jose} adapting it 
to the case of arrangements (rather than arbitrary plane curves).

\begin{dfn}\label{quasitoric} A quasi-toric 
relation of orbifold type $(3,3,3)$ is
the identity
$$F_1f^3+F_2g^3+F_3h^3=0$$
where $F_i, i=1,2,3$ is either a constant or a reduced curve such that 
$F_i=0$ are without common components and $F_1F_2F_3$ has  
the arrangement $\A$ as its zero set.
\end{dfn}
(A'priori, of course for a given arrangement there are several possibilities 
for $F_i$'s, including $F_i=1$).

\bigskip Let $E$ be elliptic curve over $\CC$ with $j$-invariant zero. 
As a model for $E$ 
we shall use the 3-fold 
ramified covering of $\CC$ branched at $0,1$ (and infinity)  
which is a closure of the affine variety in
$\CC^4$ given by the equations:
\begin{equation}\label{ellcurveformula}
   c^3=ab, a=s,b=1-s
\end{equation}
(the cyclic covering map is given by $(a,b,c,s) \rightarrow s$).
As a model of $\bar M_{\A}$ for calculation of $H^1(\bar M_{\A})$
we shall use the complete 
intersection in $\CC^5$ (with coordinates $(x,y,A,B,C)$) given by:
\begin{equation}\label{surfaceformula}
  A={{F_1(x,y)} \over {F_3(x,y)}}, B={{F_2(x,y)} \over {F_3(x,y)}},
C^3=AB
\end{equation}
(here $F_i$ are dehomogenized forms used in definition \ref{quasitoric}.
The surface (\ref{surfaceformula}) is the threefold cyclic cover of $\CC^2$
branched over the curve $F_1F_2F_3=0$ corresponding to the homomorphism 
of the fundamental group sending the classes of loops corresponding
to the curves $F_i=0$ to the same element of $\ZZ_3$. 

There is a correspondence between the $\ZZ_3$-equivariant maps 
$\phi: \bar M_{\A} \rightarrow E$ 
and quasi-toric relations given as follows. Note that such a map $\phi$
takes a fixed point of $\ZZ_3$ to one of three ramification points 
of $E \rightarrow \PP^1$. Denoting the polynomials whose zero sets 
are the remaining components of preimages 
of ramification points as $f,g,h$ respectively, one can write 
the morphism $\phi$ as 
a compactification of the restriction on affine portion of $\bar M_{\A}$
 of $\Phi: \CC^5 \rightarrow \CC^4$ given by:
\begin{equation}\label{mapformula}
a=A({f \over h})^3, b=B({g \over h})^3, c=C{{fg} \over {h^2}},
s=A({f \over h})^3
\end{equation} 

Let us define the equivalence between quasi-toric relations \ref{quasitoric}
involving $F_i,f,g,h$ and $F_i',f',g',h'$ 
saying that two such relations are equivalent if for some $\lambda, 
\lambda_f,\lambda_g,\lambda_h \in \CC$, one has
$F_i=\lambda F_i'$, ${{\lambda_f} \over {\lambda_h}} \in \mu_3,
{{\lambda_g} \over {\lambda_h}} \in \mu_3$ and 
$\lambda_f\lambda_g=\lambda_h^2$.

With this we have the following:

\begin{prop}There is one to one correspondence between 
quasi-toric relations up to above equivalence 
and non constant equivariant maps $\bar M_{\A} \rightarrow E$.
\end{prop}

\begin{proof}Indeed it follows from formulas (\ref{ellcurveformula}), 
(\ref{surfaceformula})
that $(f,g,h)$ appearing in (\ref{mapformula}) satisfy the 
quasi-toric relation. Vice versa, if $f,g,h$ satisfy a quasi-toric 
relation, the formulas (\ref{mapformula}) yield a map which 
when restricted on affine portion of $\bar M_{\A}$ provide a dominant map 
$\phi: \bar M_{\A} \rightarrow E$.
\end{proof}

As in \cite{jose} we shall view the maps $\bar M_{\A} \rightarrow E$ as 
$\CC(\bar M_{\A})$-points of the elliptic curve $E$. As such they form a group 
with the quotient by the subgroup of constant maps (the Mordell-Weil group)
being finitely 
generated (Mordell-Weil theorem cf. \cite{lang}).
The following is a special case of result shown in \cite{jose}.

\begin{theorem} 
The rank of the Mordell-Weil group $MW(\A)$ 
of $\CC(\bar M_{\A})$-points of $E$  
is equal to $2s$ where $s$ is the degree of the factor $(t^2+t+1)$ 
is (\ref{alex}). 
\end{theorem}
Finally we shall need few results on characteristic varieties of 
arrangements. Recall (cf. \cite{charvar}) that one can view 
them the components of 
affine subvariety $V$ of the torus $Char \pi_1(\PP^2-\A)$
consisting of characters $\chi$ such that 
\begin{equation}
    dim H_1(\PP^2-\A,\chi) > 0
\end{equation}
These affine subvarieties are finite unions of translated subgroups
of $Char\pi_1(\PP^2-\A)$ by the points of finite order in 
$\Char\pi_1(\PP^2-\A)$  (cf. \cite{arapura},\cite{manuscripta}).

The following results is contained in \cite{charvar} and \cite{withsergey}
respectively (cf. also \cite{dimca}):
\begin{theorem}\label{charvarpaper}
 1) The dimension of the $\omega$-eigenspace 
of the monodromy acting on $H^1(M_{\cal A})$ of 
the Milnor number $M_{\A}$ of an arrangement $\A$ is equal to 
$dim H^1(\PP^2-\A,L_{\omega})$ where $L_{\omega}$ is the local
system $\pi_1(\PP^2-\A) \rightarrow \CC^*$ sending each loop, 
having the linking number equal to one with a hyperplane in $\A$ 
and zero with all the others, 
to the 
a root of unity $\omega$.

2)The dimension of $ H^1(\PP^2-\A,L_{\chi})$
of the local system is not smaller than $d-1$ where $d$ is the dimension 
of the component of characteristic variety to which it belongs. 
\end{theorem}

More specifically, 2) follows from semi-continuity of ranks 
of cohomology of local systems (as in \cite{withsergey})
and the fact that generic 
local system in a component of positive dimension is a twist of the pull back 
of a local system on $\PP^1$ minus a collection of points. 
The latter does not change generically the rank of $H^1$ 
(cf. \cite{arapura}, Prop. 1.7) and  
relation in 2) is obvious for the local systems on complements 
to a finite union of point in $\PP^1$.

\section{Catalan Equation $A_1f^3+A_2g^3+A_3h^3=0$ 
with coefficients being products of linear forms.}
\label{catalansection}

The purpose of this section is to prove the following:

\begin{theorem}\label{catalan}
 Let $A_i=\Pi_{ j \le r_i} L_{i,j}, i=1,2,3,$ be 
products of linear forms  
$L_{i,j} \in \CC[x_1,x_2,x_3]$ ($r_i$ are non negative integers).
Assume that $L_{i,j}$ is not a multiple of $L_{i',j'}$ for any 
pair two pairs $(i,j),(i',j')$ and that lines $L_{i,j}=0$ form
an arrangement of lines with points of multiplicity 2 and 3 only.

1) The equation 

\begin{equation}\label{catformula}
A_1f^3+A_2g^3+A_3h^3=0
\end{equation}
has solution only if $A_1A_2A_3=A_1'A_2'A_3'$ where 
$A_1',A_2',A_3'$ are linear dependent over $\CC$. 
In this case the solutions of (\ref{catformula}) with linearly dependent 
$A_i$, are the pullbacks of solutions of Catalan 
equation over $\CC(t)$.

2) In the case $A_i$ are linearly dependent over $\CC$
the set of solutions to (\ref{catformula}) is infinite. 

\end{theorem}

\begin{proof} First we show that if the equation (\ref{catformula})
with $A_1=A_2=1$ has a solution then there exist factorization 
$A_3=A_1'A_2'A_3'$ such that at least two factors are non-constant 
and the Catalan equation (\ref{catformula})
with these factors as coefficients has 
a solution.
Assume that $A_1=A_2=1$, select $f,g,h$ 
which form a solution 
with  smallest possible degree of $h$ 
and such that $f,g,h$ have no common factors. 
Then $-A_3h^3=\Pi (f+\omega_i g)$
($\omega_i^3=-1$). 
Let $\tilde h$ be an irreducible factor of $h$. There is only 
one factor $f+\omega_ig$ which is divisible by $\tilde h$, since   
if there are two, then 
both $f$ and $g$ are divisible by $\tilde h$ contradicting the 
assumption of absence of common factors for $f,g,h$. 
If each irreducible factor of $h$ divides single factor  $f+\omega_ig$ 
of $f^3+g^3$, then 
$f+g=a{h'}^3A_3'$, 
$f+\omega_2g=b (h'')^3A_1',f+\omega_3g=c(h''')^3A_2'$.
So $3deg f=(deg a+deg b+deg c)+3 (deg h' +deg h''+ 3deg h''')+
deg A_3$
i.e. $a=b=c$ are constants. 
Hence $A_1'{h'}^3,A_2'{h''}^3,A_3'{h'''}^3$ 
belong to a pencil i.e. $h',h'',h'''$ are solutions to the 
Catalan equation with coefficients $A_1',A_2',A_3'$.
If $A_1',A_2'$ are constant then the degree of $h'''$ is strictly 
smaller than the degree of $h$. 
Hence such descent should lead either to Catalan equation 
with solutions having at least two non-constant coefficients or   
the Catalan equation $f^3+g^3+A_3=0$ would have a solution. 
In the latter
case, by unique 
factorization, each $f+\omega_ig$ is a product of factors 
of $A_3$ and the Catalan equation with coefficients 
$A_i''=
f+\omega_ig$ has solution in $\CC$ i.e. satisfies conclusion of 1)
since $A_i''$ are in the pencil generated by $f,g$.

Assume now that factors $A_1,A_2$ in (\ref{catformula}) 
are not constant. Let us consider the restriction of 
relation (\ref{catformula})
on one of the lines $L_{i,j}$, say $L=L_{1,j}$ which equation 
is a factor of $A_1$. 
Denoting the form obtained by restriction on $L$ of a form $M$ as $^LM$
we can write the resulting relation as
\begin{equation}\label{pencilrelation}
    ^L A_2(^Lg)^3+^L A_3(^Lh)^3=0
\end{equation}
Let $^LA_i=\Pi l_{i,j}^{a_{i,j}}$. If one $a_{i,j}>2$ then the  
zero of $l_{i,j}$ on $L$ is the point with belong to 
more than three lines $L_{i,j}=0$ forming the zero set of $A_1A_2A_3=0$
contradicting the assumption multiplicities not exceeding 3.
Hence the exponent of a factor of $^LA_i$ is $a_{i,j}=1$ or $a_{i,j}=2$. 
Due to unique factorization in $\CC[t]$ we obtain that none 
of roots of $^LA_2$ and $^LA_3$ is a  root of $g$ or $h$.
Moreover $^LA_2$ and $^LA_3$ have the same sets of zeros on $L$ 
with the same multiplicities and ${^Lg \over {^Lh}} \in \CC$.
Hence any intersection point of  $L$ and $A_i$  
also belongs to $A_{i'}$ for $i=2,3$ and $i'=2,3$.
This implies that the same is true for any point in $A_1=0$, that 
$degA_2=degA_3$ and hence the conditions of Noether $AF+BG$ 
theorem (cf. \cite{fulton}) are satisfied (the multiplicities of $A_i$ 
are equal to one). Hence 

\begin{equation}\label{pencil}
A_1=\lambda_2
A_2+\lambda_3A_3 \ \ {\rm where} \ \ {\rm deg}
 \lambda_i=0
\end{equation}

Next consider rational map $\PP^2 \rightarrow \PP^1$ given by 
$P \rightarrow (A_2(P),A_3(P))$. 
Its fibers consist of points $P$ such that 
\begin{equation}\label{ratiorelation}
{{A_2} \over {A_3}}=t \in \CC
\end{equation}
The solution (\ref{catformula}) for $A_i$ is 
also solution to
$(\lambda_2t+\lambda_3)f^3+tg^3+h^3=0$.
Since Fermat equation over the function field 
of the curve (\ref{ratiorelation}) has
only constant solutions
(cf. \cite{green}) it follows that $f,g,h$ depend only on $t$.
The same argument shows that if $a_{i,j}=2$ then
$f,g,h$ are pullbacks as well. 
Part 2) is hence proven.
\end{proof}

Next we shall consider solutions to equation (\ref{catformula}) 
over $\CC(t)$. 

\begin{theorem}\label{catalanonevar}
The solutions to (\ref{catformula}) in case $A_i \in \CC(t)$
form a group isomorphic to $Hom(Jac(C),E)$
where $C$ is the 3-fold branched cover of $\PP^1$ ramified at 
the set $\cal S$ which is the union of zeros of $A_i$ ($1 \le i \le 3$), 
$Jac(C)$ is its Jacobian 
and $E$ is the elliptic curve with an automorphism of order $3$.
In particular if $\cal S$ contains only three elements then 
the group of solutions is isomorphic to $End(E)=\ZZ^2$. 
\end{theorem}

\begin{proof} Non-constant solutions to (\ref{catformula}) correspond to 
non-constant sections of the elliptic surface (\ref{surfaceformula}).
 This elliptic surface is iso-trivial 
and is trivialized over the 3-fold 
cyclic cover $C$ of $\PP^1$ containing $t$-line as an 
open set totally ramified at the union of zeros of $A_i$ and 
possibly the infinity.
Indeed the holonomy around each zero is trivial on such 3-fold cover.
For trivial elliptic surface $C \times E$ over $C$ non constant 
sections correspond to the maps $C \rightarrow E$ by Torelli theorem 
which correspond to the maps of the Jacobian of $C$ to $E$. The last 
claim in the statement 2) and statement 3) 
follows since the Jacobian of threefold cover 
of $\PP^1$ branched at three points is $E$. 
\end{proof}

\begin{example} 1) Let 
$\omega=exp({{2 \pi \sqrt{-1}} \over 3})$ be a primitive 
root of unity of degree 3. 
A non constant solution to Catalan equation: 
$$ (t-\omega)^2(t-\omega^2)^2f^3=(1+\omega^2)(t-1)^2(t-\omega^2)^2g^3+
(-\omega^2)(t-\omega)^2
(t-1)^2h^3$$
is $$f=t-1, g=t-\omega, h=t-\omega^2$$
Indeed, after substitution we get ($1+\omega^2=-\omega$):
$$(t^3-1)^2[(t-1)-(1+\omega^2)(t-\omega)+\omega^2(t-\omega^2)]=
(t^3-1)[(t-1)+\omega(t-\omega)+\omega^2(t-\omega^2)]=0
$$

2) Suppose that $F_3=F_1+F_2$ i.e. the arrangement is composed of a pencil.
This yields the solution $f=g=h=1$ to Catalan equation.
Elliptic curve $c^3=s^2-s$ which is isomorphic to $a=s,b=1-s,c^3=ab$
is equivalent to $y^2=x^3+1$ where
$s-{1 \over 2}=y, x=4^{1 \over 3} c$.
Duplication formulas from \cite{silvermanbook} p.59, are in this case:
$$x_1:={{x^4-8x}\over {4x^3+4}}, \ \ \ \ \ 
y_1=-y-{3x^2({{x^4-8x}\over{4x^3+4}}-x)
\over {2y}}$$
(here $(x_1,y_1)$ denotes $2(x,y)$ in the sense of addition 
on elliptic curve). 
From this we infer that if we 
take $F_3=1$ (which is always the case after replacement 
$F_1 \rightarrow {{F_1}\over {F_3}}, F_2 \rightarrow {{F_2} \over {F_3}}$)
start with map $a=F_1,b=F_2$ then 
we get for doubling:
$$a:={{-F_1(F_1-2)^3} \over {(2F_1-1)^3}} \ \ \
b={{(F_1-1)(F_1+1)^3} \over {(2F_1-1)^3}}$$
i.e. 
$${f \over h}={{-F_2-1} \over {2F_1-1}} \ \ \
{g \over h}={{F_1+1} \over {2F_1-1}}$$
are solutions to Catalan equation
(indeed, one easily sees that $a+b=1$).
\end{example}

\section{On Combinatorial invariance}

\begin{proof}(Of theorems \ref{invariance} and \ref{corrected}) 
Let $\A_1$ and $\A_2$ be two 
arrangements as in the statement of the theorem. Recall that the 
characteristic polynomial of the monodromy of Milnor fiber for each  
arrangement is $(t-1)^{r(\A_i)-2}(t^2+t+1)^{s(\A_i)}$ where $r(\A_i)=
Card(\A_i)$ and $s(\A_i)$ 
is the superabundance of curves
of degree ${2d-9} \over 3$ passing through the triple points of 
the arrangement $\A_i$. Hence we need to show that $s(\A_1)>0$ 
implies  $s(\A_2)>0$ as well.

For the rest recall the results 
from \cite{jose}. Let $C$ be a possibly reducible 
plane curve for which $\omega=exp({{2 \pi \sqrt{-1}} \over 3})$
is a root of the Alexander polynomial. 
Let $\Pi_{i \in R} F_i=0$ where $R$ is the set of irreducible
components of $C$ be the reduced equation of $C$.  Then the equation 
of $C$ admits a quasi-toric relation corresponding to ordinary triple point
(cf. \cite{jose})
i.e. there exist an equality (cf. (\ref{quasitoric})): 
\begin{equation}
       \Pi_{i \in R_1} F_i f^3+\Pi_{i \in R_2} F_ig^3+
\Pi_{i \in R_3} F_ih^3=0
\end{equation}
with distinct $F_i$ and $R=\bigcup_{i=1,2,3} R_i, 
R_i \cap R_j =\emptyset, i \ne j$ is a partition 
of the set of components $R$.

Since $s(\A_1)>0$, there exist 
non trivial quasi-toric decomposition corresponding to the curve which 
is a union of lines in $\A_1$. This means that 
equations of lines of arrangement can be split into 3 groups which 
can be used to construct 
Catalan equation 
which has solutions. By theorem \ref{catalan}, $\A_1$ is composed 
of reduced pencil.

On the other hand, to each such pencil corresponds
the 2-dimensional component of characteristic 
variety of $\A_1$ containing the identity of ${\rm Char}(\pi_1(\PP^2-\A_1))$
(it consists of the pullbacks of local systems from the complement 
in $\PP^1$ to a triple of points).
Each component of characteristic variety of $\A_1$ 
containing the identity of ${\rm Char}(\pi_1(\PP^2-\A_1))$
yields a component of the resonance variety of $\A_1$ (cf. 
\cite{charvar}, \cite{withsergey}).
Recall that the resonance variety is 
 the collection of $a \in H^1(\PP^2-\A_1)$ such that the 
Aomoto complex 

\begin{equation}
    H^0(\PP^2-\A_1) \buildrel \wedge a \over \longrightarrow 
H^1(\PP^2-\A_1) \buildrel \wedge a \over \longrightarrow
H^2(\PP^2-\A_1)
\end{equation}
has non-zero cohomology at the middle term. Since such cohomology 
is determined by the cohomology algebra $H^*(\PP^2-\A_1)$ 
we obtain that the 
resonance variety of $\A_2$ has the same number of components as $\A_1$.
To each component of the resonance variety of $\A_2$ corresponds 
pencil which is reduced since the corresponding component of resonance 
variety of $\A_1$ comes from reduced pencil.
For a  map of $\PP^2-\A_2$ onto the complement in $\PP^1$ to a triple
of points $\cal P$, 
the pull back of the local system on the latter assigning to a  
standard generator of $\pi_1(\PP^1-{\cal P})$ value
 $\omega_3=exp({{2 \pi i} \over 3})$, 
is the local system on $\PP^2-\A_2$
which cohomology 
is isomorphic to the eigenspace with eigenvalue $\omega_3$ 
on $H^1(M_{\A_2})$.
Hence the claim follows. 

\end{proof}

\begin{remark}\label{dualflex}
The cases when $s>1$ correspond to arrangements which 
are composed of several pencils. One such example 
is the following.
Let $C$ be an elliptic curve and $\A$ an arrangement
of lines in plane dual to the plane containing $C$ consisting
of lines dual to nine inflection points of $C$. Then 
this arrangement has no double points, 12 triple points
and the superabundance of curves of degree $3$ 
containing these 12 points is $s=2$ (cf. \cite{dimcalehrer},\cite{nero}).
Indeed, the cohomology sequence corresponding to 
\begin{equation}
0 \rightarrow \J_S(3) \rightarrow \O_{\PP^2}(3) \rightarrow  \O_S
\rightarrow 0
\end{equation}
is the following:
\begin{equation}
 0 \rightarrow 0 \rightarrow \CC^{10} \rightarrow \CC^{12} 
\rightarrow H^1(\PP^2,\J_S(3)) \rightarrow 0
\end{equation} 
This arrangement is composed of 4 pencils (cf. \cite{charvar} section 3.3 
example 3).
For example one can use as explicit equations for this arrangement
union of three cubics 
$$x^3-y^3, x^3-z^3, y^3-z^3$$ 
or $$[(y-z)(x-z)(x-y)], [(y-\omega^2z)(z-\omega^2y)(x-\omega z)],
[(y-\omega z)(x-\omega^2 z)(x-\omega y)] \ \ \ (\omega^3=1)
$$
etc.
\end{remark}

\begin{remark}
The combinatorial invariance of $s(\A)$ would follows from the 
following refinement of the above argument describing 
the relation between $s(\A)$ and the number $l(\A)$ of pencils 
of which $\A$ is composed. This is the number of essential
components of the resonance variety of $\A$.
The number of Catalan equations corresponding to an arrangement $\A$ 
which have solutions by \ref{catalan} 
is also $l(\A)$. 
To each such pencil correspond collection the map 
$\bar M_{\cal A} \rightarrow E$ (cf. \cite{jose})
of 3-fold cyclic covers of $\PP^2$ and $\PP^1$ respectively 
defined up to automorphism of $E$ i.e.
i.e. six elements of the Mordell Weil group $MW(\A_1)$
differ by an automorphism of $E$.
Recall that $MW(\A)$ is endowed with the canonical height pairing 
(cf. \cite{serre}). 
Identification of the number of elliptic pencils on $\bar M_{\cal A}$
with the number of elements of minimal height 
in $MW(\A_1)$  and showing that the number
of elements of minimal height in $MW(\A)=\ZZ[\omega_3]^{s(\A)}$ 
is determined by $s(\A)$ 
would imply that the number of pencils which composed 
$\A$ is determined by $s(\A)$ and vice versa. Since, as follows 
from the above argument, the number of 
pencils composing the arrangement is a combinatorial invariant 
this would imply the combinatorial invariance of $s(\A)$ and hence
of the Milnor number.
\end{remark}

\end{document}